\documentclass{article}
\usepackage{graphics,amsmath,epsfig,latexsym,amssymb,cite,graphicx,a4}
\newlength{\indentedwidth} \newdimen\mathindent
\indentedwidth=\mathindent

\newtheorem{theorem}{Theorem}[section]
\newtheorem{lemma}{Lemma}[section]

\newtheorem{example}{Example}[section]

\def\ve{\varepsilon}
\def\R{\mathcal R}

\newcommand{\ov}[1]{\ensuremath{\overline{#1}}}

\DeclareMathAlphabet{\mathpzc}{OT1}{pzc}{m}{it}
\begin{document}
\vskip 0.5cm
\begin{center} {\Large \bf Representations of the quantum doubles of finite group algebras and spectral parameter dependent solutions of the Yang--Baxter equation}
\end{center}
\centerline{K.A. Dancer,%
\footnote{\tt dancer@maths.uq.edu.au} 
P.S. Isaac%
\footnote{\tt psi@maths.uq.edu.au}
and J. Links%
\footnote{\tt jrl@maths.uq.edu.au}}
~~\\
\centerline{\sl 
 Centre for Mathematical Physics, School of Physical Sciences, }
\centerline{\sl The University of Queensland, Brisbane 4072,
 Australia.} 
\vskip 0.9cm
\begin{abstract}
\vskip0.15cm
\noindent
Quantum doubles of finite group algebras form a class of quasi-triangular Hopf algebras 
which algebraically solve the Yang--Baxter equation. Each representation of the quantum double then gives a matrix solution of the Yang--Baxter equation. Such solutions do not depend on a spectral parameter, and to date there has been little investigation
into extending these solutions such that they do depend on a spectral parameter. Here we first explicitly construct the matrix elements of the generators for all irreducible representations of quantum doubles of the dihedral groups $D_n$. These results may be used to determine constant solutions of the Yang--Baxter equation. We then discuss Baxterisation ans\"atze to obtain solutions of the Yang--Baxter equation with spectral parameter and give several examples, including a new 21-vertex model. 
We also describe this approach in terms of 
minimal-dimensional representations of the quantum doubles of the alternating group $A_4$ and the symmetric group $S_4$.  
\end{abstract}

\section{Introduction}

Solutions of the Yang--Baxter equation (see (\ref{ybe}) below) provide a systematic way to construct 
exactly solvable models of two-dimensional statistical mechanics \cite{baxter}, integrable
quantum systems \cite{integrable,jm,sklyanin,lzmg}, 
as well as having applications in other areas such as knot theory \cite{knots}.
A great impetus to this field was given by Drinfeld \cite{drinfeld} who proposed the {\it quantum double
construction}. This allows for any Hopf algebra $A$ to be embedded in a larger Hopf algebra $D(A)$ 
in such a way that $D(A)$ is quasi-triangular. A consequence of the quasi-triangular property is that
there exists a canonical element $R\in D(A)\otimes D(A)$, called the {\it universal
$R$-matrix},  which solves the Yang--Baxter equation algebraically. 
Thus for any representation of $D(A)$ a matrix solution of the Yang--Baxter equation is obtained. (Below we will abuse notation and use $R$ to denote both the universal $R$-matrix and its matrix representatives.) The seminal examples
of quasi-triangular Hopf algebras were given by both Drinfeld \cite{drinfeld} and Jimbo \cite{jimbo85,jimbo86} who independently introduced
the notion of quantum algebras, which are deformations of universal enveloping algebras of Lie algebras. 

For applications to the areas mentioned above
one is generally interested in solutions of the Yang--Baxter equation with spectral parameter; i.e. for a vector 
space $V$ one looks for $ R(u,v) \in\,{\rm End}(V\otimes V) $ where $u,\,v$ are complex variables such that 
\begin{equation}
R_{12}(u,v) R_{13}(u,w) R_{23}(v,w) = 
R_{23}(v,w) R_{13}(u,w) R_{12}(u,v) \label{ybe}
\end{equation}
holds on the three-fold tensor product space $V\otimes V\otimes V$. The subscripts above refer to the 
way in which the action of
$R(u,v)$ is embedded into the space of endomorphisms on $V\otimes V\otimes V$. 
   
In the context of quantum algebras the spectral parameter arises naturally when one considers the loop representations
of affine algebras \cite{jimbo85,jimbo86}. In such instances the solutions always satisfy the {\it difference property}
$R(u,v)=R(u-v)$. However it is worth mentioning that there are solutions which do not have the difference property,
including the well-known cases of the solutions giving rise to the Hubbard model \cite{hubbard}, 
the Bariev model \cite{bariev} and the chiral Potts model
\cite{potts}. Moreover, the spectral parameter need not necessarily be a scalar, but can be a complex vector 
variable \cite{bazhanov,bandk,bkms,bgzd94,takizawa}. Below we will only concern ourselves with cases of scalar spectral parameters where the difference property does hold.

For later use we introduce the permutation operator $P$ such that $P(x\otimes y)=y\otimes x$ and set 
$\check{R}(u)=PR(u)$. Then (\ref{ybe}) can equivalently be expressed as 
\begin{equation}
\check{R}_{12}(u) \check{R}_{23}(u+v) \check{R}_{12}(v) = 
\check{R}_{23}(v) \check{R}_{12}(u+v) \check{R}_{23}(u) \label{bybe}
\end{equation}
which we will refer to as the {\it braiding} Yang--Baxter equation. It is in this form that 
the Yang--Baxter equation is relevant to knot theory \cite{knots}. Indeed setting 
$$\check{\R}=\lim_{u\rightarrow-\infty}\check{R}(u)$$
gives us 
\begin{equation}
\check{\R}_{12} \check{\R}_{23} \check{\R}_{12}= 
\check{\R}_{23} \check{\R}_{12} \check{\R}_{23} \label{braiding}
\end{equation}
which can be recognised as a defining relation in the braid group \cite{knots}. In terms of 
$${\R}=\lim_{u\rightarrow-\infty}{R}(u)$$
we have 
\begin{equation}
\R_{12} \R_{13} \R_{23} = 
\R_{23} \R_{13} \R_{12} \label{cybe}
\end{equation}
which we will refer to as the {\it constant} Yang-Baxter equation. Finally we mention that if 
\begin{equation*}
\check{R}(u)\check{R}(-u)\propto  I\otimes I 
\end{equation*} 
then the $R$-matrix is said to satisfy the {\it unitarity} condition, while if
\begin{equation}
\check{R}(0)\propto I\otimes I \label{regularity}
\end{equation} 
then it is said to satisfy the {\it regularity} condition. When the regularity condition holds there is a standard
procedure \cite{baxter,integrable,jm,sklyanin} for constructing an integrable quantum system on a one-dimensional lattice with periodic boundary conditions. 
The Hamiltonian is given by 
\begin{equation}
H = \sum_{j=1}^{L-1} h_{j,j+1} +h_{L,1}
\label{ham}
\end{equation}
where the two-site Hamiltonians are given by 
$$h=\left.\frac{d}{du} \check{R}(u)\right|_{u=0}. $$
Models constructed in this manner, and other approaches, can be solved exactly using Bethe ansatz methods 
\cite{baxter,integrable,jm,sklyanin,lzmg,knots}.

One class of quasi-triangular Hopf algebras is the set of quantum doubles of the group algebras of finite groups. 
\cite{dpr,gould93}.
(Throughout we will refer to these as {\it finite group doubles}). Applications of finite group double representations to knot theory have been addressed in \cite{tsohantjis}. These algebraic structures also underly systems of {\it anyons} in two spatial dimensions. In cases where the global symmetry of the system is spontaneously broken to a discrete gauge group, the finite group double is the appropriate structure to describe the fusion properties and statistics. The fusion properties are essentially determined by the Clebsch--Gordan decomposition of tensor products into irreducible representations. The statistics associated with the interchange of two anyons is described by braiding.  The consistency condition
for the two ways in which three anyons may be interchanged by a sequence of three two-anyon exchanges is precisely (\ref{braiding}), where
$\check{R}_{jk}$ is the operation which interchanges the $j$th and $k$th anyons. For a comprehensive review of the salient features we refer to \cite{wildbais,wild}.

Such systems may exhibit topological order \cite{toporder}, 
where quantum numbers are conserved for topological reasons, as opposed to the 
manifest symmetry. Due to the topological nature, excitations are resistant to decoherence. 
This property forms the basis of topological quantum computation which was first put forth by Kitaev 
\cite{kitaev} (see also e.g. \cite{mochon1,mochon,preskill,bonesteel}). When the symmetry is described by a finite group double, the braiding $\check{R}_{jk}$ 
is a unitary operator which can be employed as a quantum gate.     

In view of the above literature, it is surprising that there has been very little study on the role of finite group doubles in obtaining solutions
for the spectral parameter dependent Yang--Baxter equation (\ref{ybe}). Integrable systems constructed from such solutions via (\ref{ham}) 
realise models for interacting anyons with internal symmetries described by the finite group double.
Even though the models are one-dimensional there is a precedent, the Hubbard model, which leads us to believe that such models
may have applications for understanding two-dimensional systems. One property that is evident from the analysis of the Bethe ansatz solution of the one-dimemsional Hubbard model is spin-charge separation. The Hubbard model has an $so(4)\cong so(3)\oplus so(3)$ symmetry, where the two quantum numbers associated with the two copies of $so(3)$ can be assigned to spin and charge degrees of freedom. From this symmetry and the Bethe ansatz solution it can be concluded that in one dimension there exist excitations which carry spin but not charge, and vice versa, so spin-charge separation occurs \cite{kr}. It has been proposed that spin-charge separation is the mechanism responsible for high temperature superconductivity in two dimensions \cite{anderson}. Likewise, there may be insights gained into the properties of interacting anyons by studying one-dimensional models which can be solved exactly.   
  
Our aim is to investigate the extent to which solutions of the spectral parameter dependent Yang--Baxter equation can be obtained using the Hopf algebra structure of finite group doubles. This is not straightforward, as there appears to be no obvious
manner in which to consider the affine extension of a finite group double which affords loop representations. 
Using a different approach, some preliminary results in this regard have been obtained by Yang et al. \cite{kauffman}. We believe these represent the tip of an iceberg and there is ample scope for further work.
Our aim here is to continue the 
advances in this direction. 

Our starting point is to consider the quantum doubles of the dihedral groups $D_n$. Of all non-abelian finite groups, the series of dihedral groups has the simplest representation theory. We will show that the representation theory for the quantum doubles is also readily tractable. Using the general results on the representation theory of finite group doubles given in \cite{dpr,gould93}, we begin by explicitly constructing {\it all} irreducible representations for the doubles $D(D_n)$\footnote{After completing this work we learned that for odd $n$ these representations have been constructed in the thesis by de Wild Propitius\cite{wild}. Results for even $n$ have independently been obtained by Slingerland \cite{joost}} .
{}From these results it is straightforward to explicitly construct solutions  $R$ for the constant Yang--Baxter equation
(\ref{cybe}) which do not depend on the spectral parameter. 

Our next goal is to determine if these constant solutions of the Yang--Baxter equation can be extended to 
spectral parameter dependent solutions. This is a procedure colloquially known as {\it Baxterisation} as coined by Jones \cite{jones}, and there is a sizable literature on this topic \cite{cgx,zgb,dgz94,dgz96,grimm,maillard,gm}. We begin by studying the case of the two-dimensional 
irreducible representations of $D(D_n)$ and find that Baxterisation can be performed successfully. Our approach is based on an ansatz taken from \cite{cgx} which is chosen by symmetry considerations. In all these cases we find that 
the resulting solution of the Yang-Baxter equation is a particular case of the well-known trigonometric six-vertex model
in the symmetric gauge at a specific root of unity. We metion that this result is not obvious in the sense that 
the Baxterisation is not underpinned by a Hecke algebra representation.  

We then turn our attention to the three-dimensional irreducible representations of $D(D_n)$. We find that the only cases for which an irreducible three-dimensional representation exists are $D(D_3)$ and $D(D_6)$. All instances give unitarily equivalent 
constant solutions of the Yang--Baxter equation. Our Baxterisation ansatz leads to a 21-vertex solution of the spectral parameter dependent Yang--Baxter equation, which as far as we can ascertain is new.   

Rather than continuing on to investigate higher dimensional representations of the $D(D_n)$ series, we finish by considering
minimal-dimensional representations of the double of the alternating group $A_4$ and the symmetric group $S_4$. Neither of these 
cases admit irreducible two-dimensional representations, but they both admit three-dimensional ones. Our ansatz 
for Baxterising the constant solutions does lead us to spectral parameter dependent solutions.
Remarkably though, we find in these latter examples that the infinite spectral parameter limit of the Baxterised solutions do not necessarily reproduce the original solutions of the constant Yang--Baxter equation.

\section{The dihedral group $D_n$}

Consider the dihedral group $D_n$.  This has two generators $\sigma, \tau$ satisfying:

$$ \sigma^n = e,\; \tau^2 = e,\; \tau \sigma = \sigma^{n-1} \tau .$$

\noindent The properties of $D_n$ vary according to whether $n$ is odd or even, with the odd case being slightly simpler.

\subsection{$D_n$ where $n$ is odd}

When $n$ is odd, there are $({n+3})/{2}$ conjugacy classes divided into three families, given by:

\begin{align*}
&\{ e \}, \\
&\{ \sigma^k, \sigma^{-k} \} \quad \text{for } 1 \leq k \leq \frac{n-1}{2}, \\
&\{ \sigma^i \tau, \, 0 \leq i \leq n-1 \}.\notag
\end{align*}
There are $({n+3})/{2}$ irreducible representations (irreps), 
two of which are one-dimensional and the remaining $({n-1})/{2}$ which are two-dimensional. They are given by:

$$ \pi_\pm (\sigma) =1,\quad \pi_\pm (\tau) = \pm 1  $$ 

\noindent and

\begin{equation*}
\pi_k (\sigma) = \begin{bmatrix}
                    \omega^k & 0  \\ 0 & \omega^{-k}
                 \end{bmatrix}, \quad
\pi_k (\tau) = \begin{bmatrix}
                 0 & 1 \\ 1 & 0 
               \end{bmatrix}, \quad
\omega = \exp\left({\frac{2 \pi i}{n}}\right), \; 1 \leq k \leq \frac{n-1}{2}.
\end{equation*}

\noindent As required, the sum of the squares of the dimensions of the irreps 
is $2 \times 1^2 + (n-1)/2 \times 2^2 = 2n = |D_n|$ \cite{curtis}.

\subsection{$D_n$ where $n$ is even}

When $n$ is even, there are $({n}+6)/{2}$ conjugacy classes divided into five families, given by:

\begin{align*}
&\{e\}, \notag \\
&\{\sigma^{{n}/{2}} \}, \notag \\
&\{ \sigma^k, \sigma^{-k} \} \quad \text{for } 1 \leq k \leq {(n-2)}/{2}, \\
&\{ \sigma^{2j} \tau, \; 0 \leq j \leq {(n-2)}/{2} \}, \notag \\
&\{ \sigma^{(2j+1)} \tau, \; 0 \leq j \leq {(n-2)}/{2} \}. \notag 
\end{align*}

\noindent The $({n+6})/{2}$ irreps consist of 4 one-dimensional irreps and ${(n-2)}/{2}$ two-dimensional irreps.  They are given by:

$$ \pi(\sigma) = (-1)^a, \quad \pi (\tau) = (-1)^b \quad \text{for } \; a,b \in \{ 0,1 \} $$

\noindent and 

\begin{equation*}
\pi_k (\sigma) = \begin{bmatrix}
                    \omega^k & 0  \\ 0 & \omega^{-k}
                 \end{bmatrix}, \quad
\pi_k (\tau) = \begin{bmatrix}
                 0 & 1 \\ 1 & 0 
               \end{bmatrix}, \quad
\omega = \exp\left({\frac{2\pi i}{n}}\right),\; 1 \leq k \leq \frac{n-2}{2}.
\end{equation*}

\noindent Again, the sum of the squares of the dimensions of the irreps is $2n$.


\section{The quantum double algebra $D(G)$}

\noindent Here we give a brief survey of finite group doubles $D(G)$. Throughout, our approach follows that of 
\cite{gould93}. 
Let $A$ be the group algebra of a finite group $G$ over the complex field $\mathbb{C}$.  Then $A$ becomes a co-commutative Hopf algebra with coproduct, antipode and counit respectively defined by

$$ \Delta(g) = g \otimes g, \quad S(g) = g^{-1}, \quad \varepsilon(g) =e, \quad \forall g \in G.$$

\noindent Let $A^*$ be the dual space of $A$, so $A^* = \{f|f:A \rightarrow \mathbb{C}\}$.  Then $A^*$ becomes an algebra on the dual elements $g^*$ defined by 

$$ g^*(h) = \delta(g,h) \quad \forall g,h \in G,$$ 

\noindent which have the property

$$g^* h^* = \delta(g,h) h^*.$$
 The resultant dual algebra is commutative and does not have an interesting representation theory.  Now we follow the quantum double construction to obtain $D(G)$, which is a $|G|^2$-dimensional algebra spanned by the free products

$$ g h^*, \quad g, h \in G.$$

\noindent The elements $h^* g$ are calculated using

$$ h^* g = g (g^{-1}hg)^*.$$
Then $D(G)$ becomes a quasi-triangular Hopf algebra with coproduct $\ov{\Delta}$, antipode $\ov{S}$ and counit $\ov{\varepsilon}$ given by:

\begin{align*}
&\ov{\Delta}(gh^*) = \sum_{k \in G} g (k^{-1}h)^* \otimes gk^* = \sum_{k \in G}
  gk^* \otimes g(kh^{-1})^*, \\
&\ov{S}(gh^*) = g^{-1} (g h^{-1} g^{-1})^*, \\
&\ov{\ve}(gh^*) = \delta(h,e).
\end{align*}
Note that we identify $g \ve$ with $g$ and $e g^*$ with $g^*$ for all $g \in G$.
The universal $R$-matrix is given by 

$$ R = \sum_{g \in G} g \otimes g^* $$

\noindent which can easily be shown to satisfy the defining relations for a quasi-triangular Hopf algebra:
\begin{eqnarray}
R\ov{\Delta}(a)&=&\ov{\Delta}^T(a)R, \quad \forall\,a\in D(G), \label{qt1}\\
(\ov{\Delta}\otimes {\rm id}) R&=&R_{13}R_{23},   \label{qt2}\\
({\rm id}\otimes \ov{\Delta})R&=& R_{13}R_{12}, \label{qt3}
\end{eqnarray}  
where $\ov{\Delta}^T$ is the opposite coproduct
$$ \ov{\Delta}^T(gh^*) = \sum_{k \in G} gk^* \otimes g (k^{-1}h)^* = \sum_{k \in G}
g(kh^{-1})^* \otimes   gk^*. $$ 
It follows from the relations (\ref{qt1},\ref{qt2},\ref{qt3}) that $R$ is a solution of the constant 
Yang-Baxter equation.
We note that in any tensor product representation $\pi\otimes \pi$ we have that $\check{R}=PR$ commutes 
with the action of the finite group double; i.e. 
$$[\check{R},\,(\pi\otimes \pi) \ov{\Delta}(a)]=0 \quad \forall \,a\in D(G). $$ 

\section{Representation theory of $D(G)$}

Several properties of group algebras extend to the quantum double.  For example, the set 

$$Q = \{ gh^* | g, h \in G \text{ with } gh=hg \},$$
is stable under the adjoint action of $G$, i.e. $gQg^{-1}=Q$.  Hence $Q$ can be partitioned into $G$-conjugacy classes, which implies \cite{gould93}

\begin{theorem}
The number of non-isomorphic $D(G)$-modules equals the number of $G$-equivalence classes of $Q$.
\end{theorem}
Moreover, a construction for these modules is known \cite{gould93}.  The general form is included in this section, with the explicit results for the odd and even dihedral groups given in the following two sections respectively.  

First partition $G$ into conjugacy classes

$$ G = \bigcup_k C_k.$$

\noindent Recall that the centraliser subgroup of an element $h$ is defined by

$$ Z(h) = \{ g \in G | gh =hg\}.$$
Then for each conjugacy class $C_k$ choose a representative $g_k \in C_k$ and set $Z_k = Z(g_k)$ to be the centraliser subgroup of $g_k$, noting that $|Z_k| |C_k| = |G|$.  Denote the group algebra of $Z_k$ by $A_k$.  
Also, for each $s \in C_k$ choose a fixed element $\alpha_s \in G$ satisfying

$$s = \alpha_s g_k \alpha_s^{-1}.$$

\noindent For simplicity choose $\alpha_{g_k} = e \; \forall g_k$. 

\begin{lemma} \label{alphas}
We have the following properties of $\alpha_s$: 
\begin{enumerate}
\item$ G = \bigcup_{s \in C_k} \alpha_s Z_k$
\item Given $g \in G, s \in C_k, \exists t \in C_k$ unique with the property
$ \alpha_t^{-1} g \alpha_s \in Z_k;$
explicitly $t = gsg^{-1}$.
\end{enumerate}
\end{lemma}

\noindent Again, a proof can be found in \cite{gould93}.  

The irreducible modules of $D(G)$ can be constructed from modules of the group algebras $A_k$.  Let $V^k_\beta$ denote an irreducible $A_k$ module.  Then there is a corresponding induced $A$-module \cite{curtis}

$$V_{k,\beta} \subseteq A \otimes_{A_k} V^k_\beta$$

\noindent spanned by vectors

$$ v(s) = \alpha_s \otimes v, \qquad v \in V^k_\beta, s \in C_k$$

\noindent where the action of $G$ is given by

$$g (\alpha_s \otimes v) = \alpha_{gsg^{-1}} \otimes ( \alpha_{gsg^{-1}}^{-1}
  g \alpha_s)v$$

\noindent or, equivalently, 

$$gv(s) = \bigl( \alpha_{gsg^{-1}}^{-1}g \alpha_s v \bigr) (gsg^{-1}).$$
Note ${\rm dim} V_{k,\beta} = |C_k|\, {\rm dim} V^k_{\beta}$.  
It follows from Lemma \ref{alphas} that $V_{k, \beta}$ is an $A$-module under this definition.

The module $V_{k,\beta}$ can be decomposed according to

$$V_{k,\beta} = \bigoplus_{s \in C_k} V_{k, \beta} (s)$$ 

\noindent where 

$$V_{k, \beta}(s) = \{ v(s) | v \in V^k_\beta \}.$$

\noindent The latter becomes an irreducible module over the group algebra of $Z(s) = \alpha_s Z_k \alpha_s^{-1}$.  When $s=g_k$ the module is isomorphic to $V^k_\alpha$. Then 
$V_{k, \beta}$ becomes an irreducible $D(G)$ module with the action:

$$ h^* v(s) = \delta(h,s) v(s), \qquad \forall h \in G.$$

\noindent Moreover, two $D(G)$ modules of this form, $V_{k, \beta}, V_{l, \gamma}$, are isomorphic iff $k=l$ and $V^k_\beta, V^k_\gamma$ are isomorphic.  Then using counting arguments it can be shown that \cite{gould93}:

\begin{theorem}
Every irreducible $D(G)$-module is isomorphic to one of the $V_{k, \beta}$.
\end{theorem}



\section{Representations of $D(G)$ where $G=D_n$, $n$ even} \label{irrepseven}

The conjugacy classes $C_k$ of $G=D_n$, chosen representatives $g_k$, corresponding centraliser subgroups 
and the elements $\alpha_s,\, \forall s \in C_k,$ are given below in Table \ref{even}.

\begin{table}[ht]
\caption{$C_k, g_k$, $Z_k$ and $\alpha_s$ for $G=D_n$, $n$ even.} \label{even}
\centering
\begin{tabular}{|l|l|l|l|} \hline
$C_k$ & $g_k$ & $Z_k = Z(g_k)$ & $\alpha_s, \forall s \in C_k$ \\ \hline 
$\{e\}$ & $e$ & $D_n$ & $\alpha_e=e$\\
$\{ \sigma^{\frac{n}{2}} \}$ & $\sigma^{\frac{n}{2}}$ & $D_n$ & 
  $\alpha_{\sigma^{n/2}} = e$ \\
$\{\sigma^k, \sigma^{-k} \}, 1 \leq k < \frac{n}{2}$ & $\sigma^k$ & 
  $\{\sigma^i \,|\, 0 \leq i < n \}$ & $\alpha_{\sigma^k}=e, 
  \alpha_{\sigma^{-k}} = \tau$ \\
$\{ \sigma^{2i} \tau \,|\, 0 \leq i < \frac{n}{2} \}$ & $\tau$ & 
  $\{ e, \tau, \sigma^{\frac{n}{2}}, \sigma^{\frac{n}{2}}\tau \}$ & 
  $\alpha_{(\sigma^{2i} \tau)} = \sigma^{i}$\\ 
$\{ \sigma^{2i+1} \tau \,|\, 0 \leq i < \frac{n}{2} \}$ & $\sigma \tau$ & 
  $\{ e, \sigma \tau, \sigma^{\frac{n}{2}}, \sigma^{\frac{n}{2}+1}\tau \}$ & 
  $\alpha_{(\sigma^{2i+1} \tau)} = \sigma^{i}$\\ \hline
\end{tabular}
\end{table}

Throughout the remainder of this paper $E^i_j$ denotes an elementary matrix with a 1 in 
the $(i,j)$ position and zeroes elsewhere.  We also abuse notation by using $g$ to denote both an element of the algebra $D(G)$ and its matrix representative in a given irrep, which should be clear from the context.

\

\noindent \underline{Representations induced by $C_k = \{e\}$}

\

The module elements are of the form $e \otimes v$ where $v \in V$, $V$ a $D_n$-module.  In representation terms, there are 4 one-dimensional irreps and $(\frac{n}{2}-1)$ two-dimensional irreps.  They are:

$$\sigma=(-1)^a, \quad \tau = (-1)^b, \quad g^* = \delta(g, e)$$

\noindent where $a, b \in \{0,1\}$, and 

\begin{equation*}
 \sigma = \begin{bmatrix}
                    \omega^k & 0  \\ 0 & \omega^{-k}
                 \end{bmatrix}, \quad
\tau = \begin{bmatrix}
                 0 & 1 \\ 1 & 0 
               \end{bmatrix}, \quad 
g^* = \delta(g, e) I_2 
\end{equation*} 

\noindent where $1 \leq k < \frac{n}{2}$.

\noindent \underline{Representations induced by $C_k = \{\sigma^{{n}/{2}}\}$}

\

The module elements are of the form $ e \otimes v$ where $v \in V$, $V$ a $D_n$-module.  Again, there are 4 one-dimensional irreps and $({n-2})/{2}$ two-dimensional irreps.  They are:

$$\sigma=(-1)^a, \quad \tau = (-1)^b, \quad g^* = 
  \delta(g, \sigma^{{n}/{2}})$$

\noindent where $a, b \in \{0,1\}$, and 

\begin{equation*}
 \sigma = \begin{bmatrix}
                    \omega^k & 0  \\ 0 & \omega^{-k}
                 \end{bmatrix}, \quad
\tau = \begin{bmatrix}
                 0 & 1 \\ 1 & 0 
               \end{bmatrix}, \quad 
g^* = \delta(g, \sigma^{{n}/{2}}) I_2 
\end{equation*} 

\noindent where $1 \leq k < {n}/{2}$.

\

\noindent \underline{Representations induced by $C_k = \{\sigma^k, \sigma^{-k}\}, \; 1 \leq k < {n}/{2}$}

\

The module elements are of the form $e \otimes v, \tau \otimes v$ where $v \in V$, $V$ a module of the group algebra of $Z_k=\{ \sigma^j\, | \, 0 \leq j < n \}$.  There are $n$ such $A_k$-modules, with the corresponding representations given by $\pi (\sigma) = \omega^j,\, 0 \leq j < n$ where $\omega = \exp({{2\pi i}/{n}})$.  Thus we have $n({n}-2)/2$ different irreducible representations of $D(D_n)$ induced by these conjugacy classes, given by:

\begin{equation*}
 \sigma = \begin{bmatrix}
                    \omega^j & 0  \\ 0 & \omega^{-j}
                 \end{bmatrix}, \quad
\tau = \begin{bmatrix}
                 0 & 1 \\ 1 & 0 
               \end{bmatrix}, \quad 
(\sigma^k)^* = \begin{bmatrix}
                 1 & 0 \\ 0 & 0
               \end{bmatrix}, \quad
(\sigma^{-k})^* = \begin{bmatrix}
                 0 & 0 \\ 0 & 1
               \end{bmatrix}, \quad
g^* = 0 \text{ otherwise},
\end{equation*} 

\noindent where $0 \leq j < n, \, 1 \leq k <{n}/{2}$.

\

\noindent \underline{Representations induced by $C_k = \{\sigma^{2j} \tau \,|\, 0 \leq j  < {n}/{2} \}$}

The module elements are of the form $\sigma^{j} \otimes v,\: 0 \leq j < 
 {n}/{2}$, where $v \in V$, $V$ a module of the group algebra of 
 $Z_k = \{ e, \tau, \sigma^{{n}/{2}}, \sigma^{{n}/{2}}\tau  \}$.  
 Hence there are four $(n/2)$-dimensional irreps of this form.  They are:

\begin{align*}
&\sigma = A \in M_{\frac{n}{2} \times \frac{n}{2}} \,\,\,\quad\quad\quad\quad \text{ where } 
  [A]_{ij} = (-1)^{a\, \delta (i, 1)} \delta(i, j+1), \quad  
  \text{addition } mod \; n/2 \\
&\tau = (-1)^{a+b} B \in M_{\frac{n}{2} \times \frac{n}{2}} \quad 
  \text{ where } [B]_{ij} = (-1)^{a \delta (i,1)} \delta(i+j, 2), 
  \quad \text{addition } mod \; n/2 \\
&(\sigma^i)^* = 0, \quad (\sigma^{2j} \tau)^* = E^{j+1}_{j+1}, 
  \quad (\sigma^{(2j+1)}\tau)^* = 0, \quad 0 \leq i < n,\; 0 \leq j < n/2
\end{align*}

\noindent where $a, b \in \{ 0,1 \}$.

\begin{example}
In $D(D_6)$, $\sigma$ and $\tau$ are as follows:
\begin{equation*}
\sigma = \begin{bmatrix}
      0 & 0 & (-1)^a \\ 1 & 0 & 0 \\0 & 1 & 0  
          \end{bmatrix}, \quad
\tau = (-1)^{b} \begin{bmatrix}
      1 &0& 0 \\ 0 & 0 & (-1)^{a} \\ 0 &  (-1)^{a} & 0
       \end{bmatrix}
\end{equation*}
\noindent where $a,b \in \{0,1 \}$, giving 4 three-dimensional irreps.
\end{example}

\noindent \underline{Representations induced by $C_k = \{\sigma^{2j+1} \tau \,|\, 0 \leq j  < {n}/{2} \}$}

The module elements are of the form $\sigma^{j}\, \otimes v, \: 0 \leq j < 
 {n}/{2}$, where $v \in V$, $V$ a module of the group algebra of 
 $Z_k = \{ e, \sigma \tau, \sigma^{{n}/{2}}, \sigma^{{(n+2)}/{2}}\tau\}$.
 Hence there are four $(n/2)$-dimensional irreps of this form.  They are:

\begin{align*}
&\sigma = A \in M_{\frac{n}{2} \times \frac{n}{2}} \quad\,\,\,\quad\quad\quad \text{ where } 
  [A]_{ij} = (-1)^{a\, \delta (i, 1)} \delta(i, j+1), \quad\,\,  \text{addition } 
  mod \; n/2 \\
&\tau = (-1)^{a+b} B \in M_{\frac{n}{2} \times \frac{n}{2}} \quad 
  \text{ where } [B]_{ij} = \delta(i+j, n/2+1), \quad\quad\quad \text{addition } 
  mod \; n/2 \\
&(\sigma^i)^* = 0, \quad (\sigma^{2j} \tau)^* = 0, \quad (\sigma^{(2j+1)} 
  \tau)^* = E^{j+1}_{j+1}, \quad 0 \leq i < n,\; 0 \leq j < n/2
\end{align*}

\noindent where $a, b \in \{ 0,1 \}$.

\begin{example}
In $D(D_6)$, $\sigma$ and $\tau$ are as follows:
\begin{equation*}
\sigma = \begin{bmatrix}
      0 & 0 & (-1)^a \\ 1 & 0 & 0 \\0 & 1 & 0  
          \end{bmatrix}, \quad
\tau = (-1)^{a+b} \begin{bmatrix}
      0 &0& 1 \\ 0 & 1 & 0 \\ 1 & 0 & 0
       \end{bmatrix}
\end{equation*}
\noindent where $a,b \in \{0,1 \}$, giving 4 three-dimensional irreps.
\end{example}

Thus when $n$ is even $D(D_n)$ has 8 one-dimensional irreps, 8 $({n}/{2})$-dimensional irreps and 
\mbox{$(n+2)(n/2-1)$} two-dimensional irreps.  The sum of the squares of the dimension of the irreps is $|D_n|^2 = |D(D_n)|$, as required \cite{gould93}. 

Analogous results for odd $n$ are given in the Appendix.


\section{Solutions of the Yang--Baxter equation associated with the two-dimensional irreps of $D(D_n)$}

As stated earlier, a solution to the constant Yang--Baxter equation in $D(G)$ is given by

$$ R = \sum_{g \in G} g \otimes g^*.$$
In this section we will only consider the two-dimensional irreps of $D(G)$.  Then, by inspection, the $R$-matrix will always be of the form 

$$R = {\rm diag} (\omega^j, \omega^{-j},\omega^{-j}, \omega^j )$$

\noindent where $\omega = \exp({2\pi i}/{n})$ and $0 \leq j < n$, with the cases $j = 0$ and $j = {n}/{2}$ being trivial.  Thus $\check{R} = PR$ is given by:

\begin{equation} 
\check{R} = \begin{bmatrix}
             \omega^j & 0 & 0 &0 \\ 
	     0 & 0 & \omega^{-j} & 0 \\
	     0 & \omega^{-j} & 0 & 0 \\
	     0 & 0 & 0 & \omega^j
          \end{bmatrix}. \label{2d} 
\end{equation}

As remarked earlier, $\check{R}$ commutes with the action of $D(D_n)$. We seek a solution $\check{R}(x)$ of the braiding Yang--Baxter equation which also has this symmetry. Now the matrix (\ref{2d}) has 3 different eigenvalues, namely $w^j, \pm w^{-j}$.  Hence for any $k$, $\check{R}^k$ can be written as a linear combination of $I \otimes I, \check{R}$ and $\check{R}^{-1}$. Making a change of variables  $x=\exp(u),\,z=\exp(v)$ in equation \eqref{bybe}, we look for $\check{R}(x)$ in the following form:

$$ \check{R}(x) = f(x) I \otimes I + g(x) \check{R} + h(x) \check{R}^{-1}.$$

\noindent Then 

\begin{equation*}
\check{R}(x) = \begin{bmatrix}
             A(x) & 0 & 0 &0 \\ 
	     0 & f(x) & B(x) & 0 \\
	     0 & B(x) & f(x) & 0 \\
	     0 & 0 & 0 & A(x)
          \end{bmatrix}
\end{equation*}

\noindent where 

\begin{align*}
&A(x) = f(x) + \omega^j  g(x) + \omega^{-j} h(x), \\
&B(x) = \omega^{-j} g(x) + \omega^j h(x).
\end{align*}
We directly apply the braiding Yang--Baxter equation in the form

\begin{equation}
 \check{R}_{12}(x) \check{R}_{23}(xz) \check{R}_{12}(z) 
= \check{R}_{23}(z) \check{R}_{12}(xz) \check{R}_{23}(x). \label{mbybe}
\end{equation}

\noindent Although there are 20 non-zero entries on each side of the equation, there are only 2 independent non-trivial relations that must be satisfied.  These are:

\begin{align}
A(z) f(xz) A(x) = f(x) A(xz) f(z) + B(x) f(xz) B(z), \label{AfA} \\
A(z) B(xz) f(x) = f(x) A(xz) B(z) + B(x) f(xz) f(z). \label{ABf}
\end{align}

\noindent Note that $f(x) = 0 \; \forall x$ trivially satisfies the Yang--Baxter equation. 

A proposal for constructing $\check{R}(x)$ when $\check{R}$ has three distinct eigenvalues $\lambda_1, \lambda_2$ and $\lambda_3$ has been discussed in \cite{cgx,kauffman}, but it has not been proven to always be true.  The conjecture is 

$$\check{R}(x) = (\lambda_1 + \lambda_2 + \lambda_3 + \lambda_1 \lambda_3 
 \lambda_2^{-1}) x I \otimes I - (x-1)\check{R} + \lambda_1 \lambda_3 x(x-1) 
 \check{R}^{-1}. $$

Three distinct solutions are obtained by changing the ordering of the eigenvalues. We note that when the ansatz holds we
have 
\begin{equation*}
\check{\R}=\lim_{x\rightarrow 0} R(x)= R. 
\end{equation*}

\noindent Applying this ansatz to (\ref{2d}) we find that if $\lambda_2 = \pm \omega^{-j}$ then $f(x) = 0$, which we have already shown gives a trivial result.  Hence we consider the case when $\lambda_2 = \omega^j$.  This gives:

\begin{align*}
&f(x) = (\omega^j - \omega^{-3j})x, \quad g(x) = -(x-1), \quad h(x) = 
  -w^{-2j} x(x-1) \\
\Rightarrow \quad & A(x) = \omega^j - \omega^{-3j} x^2, \quad 
 B(x) = -\omega^{-j}  (x^2-1) .
\end{align*}

It can be easily shown that $f(x), A(x)$ and $B(x)$ satisfy relations \eqref{AfA} and \eqref{ABf}.  Hence we have solutions to the braiding Yang--Baxter equation, which are 

\begin{equation*}
\check{R}(x) = \begin{bmatrix}
                  \omega^j - \omega^{-3j} x^2 & 0 & 0 & 0 \\
		  0 & (\omega^j - \omega^{-3j})x & -\omega^{-j} (x^2-1) & 0 \\
		  0 & -\omega^{-j} (x^2-1) & (\omega^j - \omega^{-3j})x & 0 \\
		  0 & 0 & 0 & \omega^j - \omega^{-3j}x^2
		\end{bmatrix}
\end{equation*}

\noindent where $\omega = \exp({{2\pi i}/{n}})$ and $0 \leq j < n$.  Rescaling by a factor of $\omega^j x^{-1}$, we can write:

\begin{equation}
\check{R}(x) = \begin{bmatrix}
                  \omega^{2j}x^{-1} - \omega^{-2j} x & 0 & 0 & 0 \\
		  0 & \omega^{2j} - \omega^{-2j} & x^{-1}-x & 0 \\
		  0 & x^{-1}-x & \omega^{2j} - \omega^{-2j} & 0 \\
		  0 & 0 & 0 & \omega^{2j}x^{-1} - \omega^{-2j} x
		\end{bmatrix}. \label{6v}
\end{equation}

\noindent Note that the unitarity condition $\check{R}(x) \check{R}(x^{-1}) 
= [\omega^{4j}+\omega^{-4j} -(x^2 + x^{-2})] I \otimes I$ is satisfied. We can recognise (\ref{6v}) as
specific cases of
the six-vertex solution in the symmetric gauge, where the parameter $q$ in the general solution 
is constrained to be a root of unity. We remark that the choice of gauge is related to the gradation chosen for the 
affine algebra (e.g. see \cite{gz}). It is interesting to note that in the non-symmetric gauge 
the constant solution
$\check{\R}=\lim_{x\rightarrow 0} R(x)$ can be used to give rise to a representation of the Temperley-Lieb algebra. In this case there is a well known procedure for Baxterising $\check{\R}$ to recover $\check{R}(x)$ \cite{jones,bk}. This is not the case for the symmetric gauge case described above.  

\section{Solutions of the Yang--Baxter Equation associated with the three-dimensional irreps of $D(D_3)$ and $D(D_6)$}

{}From the construction given in Sections \ref{irrepsodd}, \ref{irrepseven} above, we find that three-dimensional 
irreps only occur for $D(D_3)$ and $D(D_6)$. Moreover, the two three-dimensional irreps of $D(D_3)$ are also representations for the $D(D_3)$-subalgebra of $D(D_6)$, so these cases give identical $R$-matrices.   
For any three-dimensional irrep of $D(D_6)$ we find $\check{R}=PR=P\sum_g g \otimes g^*$ is of the following form:

\begin{equation*}
\check{R}=(-1)^b \begin{bmatrix}
  1 & 0 & 0 & 0 & 0 & 0 & 0 & 0 & 0 \\
  0 & 0 & 0 & 0 & 0 & 1 & 0 & 0 & 0 \\
  0 & 0 & 0 & 0 & 0 & 0 & 0 & (-1)^a & 0 \\
  0 & 0 & (-1)^a & 0 & 0 & 0 & 0 & 0 & 0 \\
  0 & 0 & 0 & 0 & 1 & 0 & 0 & 0 & 0 \\
  0 & 0 & 0 & 0 & 0 & 0 & (-1)^a & 0 & 0 \\
  0 & (-1)^a & 0 & 0 & 0 & 0 & 0 & 0 & 0 \\
  0 & 0 & 0 & 1 & 0 & 0 & 0 & 0 & 0 \\
  0 & 0 & 0 & 0 & 0 & 0 & 0 & 0  & 1
      \end{bmatrix}
\end{equation*}

\noindent where $a,b \in \{0,1\}$.  Without loss of generality we take $b=0$, and find the eigenvalues of $\check{R}$ are $1, \omega$ and $\omega^2$ with multiplicities $5,2$ and $2$ respectively, where $\omega = \exp({2\pi i}/{3})$. As the two values of $a$ give unitarily equivalent $R$-matrices we choose to take $a=0$.  We can then write 

$$\check{R}(x) = f(x) I \otimes I + g(x) \check{R} + h(x) \check{R}^{-1}.$$

\noindent Again we follow the procedure outlined in \cite{cgx,kauffman} to find possible solutions, which gives the following:

\begin{table}[ht]
\caption{Possible solutions for $f(x), g(x)$ and $h(x)$.}
\centering
\begin{tabular}{|c|c|c|} \hline
$f(x)$ & $g(x)$ & $h(x)$ \\ \hline
$x$& $1-x$& $x(x-1)$  \\
$\omega x$& $1-x$& $\omega^2 x(x-1)$  \\
$\omega^2 x$ & $1-x$ & $\omega x(x-1)$ \\ \hline
\end{tabular}
\end{table}

Using Mathematica we find only the first of these possible solutions satisfies the braiding Yang--Baxter equation
(\ref{mbybe}):

\noindent 
{\small
\begin{equation*}
\check{R}(x)=\begin{bmatrix}
  x^2-x+1 & 0 & 0 & 0 & 0 & 0 & 0 & 0 & 0 \\
  0 & x & 0 & 0 & 0 & 1-x &  x(x-1) & 0 & 0 \\
  0 & 0 & x &  x(x-1) & 0 & 0 & 0 &  1-x & 0 \\
  0 & 0 &  1-x & x & 0 & 0 & 0 & x(x-1) & 0 \\
  0 & 0 & 0 & 0 & x^2-x+1 & 0 & 0 & 0 & 0 \\
  0 & x(x-1) & 0 & 0 & 0 & x & 1-x & 0 & 0 \\
  0 & 1-x & 0 & 0 & 0 &  x(x-1) & x & 0 & 0 \\
  0 & 0 &  x(x-1) & 1-x & 0 & 0 & 0 & x & 0 \\
  0 & 0 & 0 & 0 & 0 & 0 & 0 & 0  & x^2-x+1
      \end{bmatrix}.
\end{equation*}
}

\noindent The corresponding $R$-matrix $R(x)= P\check{R}(x)$ is given by

\noindent 
{\small
\begin{equation}
R(x)=\begin{bmatrix}
  x^2-x+1 & 0 & 0 & 0 & 0 & 0 & 0 & 0 & 0 \\
  0 & 0 & 1-x & x & 0 & 0 & 0 & x(x-1) & 0 \\
  0 & 1-x & 0 & 0 & 0 & x(x-1) & x & 0 & 0 \\
  0 & x & 0 & 0 & 0 & 1-x & x(x-1) & 0 & 0 \\
  0 & 0 & 0 & 0 & x^2-x+1 & 0 & 0 & 0 & 0 \\
  0 & 0 & x(x-1) & 1-x & 0 & 0 & 0 & x & 0 \\
  0 & 0 & x & x(x-1) & 0 & 0 & 0 & 1-x & 0 \\
  0 & x(x-1) & 0 & 0 & 0 & x & 1-x & 0 & 0 \\
  0 & 0 & 0 & 0 & 0 & 0 & 0 & 0  & x^2-x+1
      \end{bmatrix}.
\label{newR}
\end{equation}
}

\noindent  Note $\check{R}(x)\check{R}(x^{-1}) = (x-1+1/x)^2 \, I \otimes I$, so the unitarity property holds.  

The above solution gives rise to a 21-vertex model, which appears to be new. It does not belong to the class of 21-vertex models discussed in \cite{golzer}. Viewed as a two-dimensional lattice 
statistical mechanics model though, it does not have real, non-negative Boltzmann weights. Since
the regularity property (\ref{regularity}) holds, we can however construct an integrable one-dimensional model. Even though $\check{R}(u)$ is not Hermitian, we obtain a Hermitian Hamiltonian
in the following manner.  We rescale $\check{R}(u)$ by a factor of $i/x$ and define the two-site Hamiltonian $h$ as

\begin{equation*}
h = \frac{d}{dx}\left.\frac{i\check{R}(x)}{x}\right|_{x=1} =  
    i \begin{bmatrix}
       0 & 0 & 0 & 0 & 0 & 0 & 0 & 0 & 0 \\
       0 & 0 & 0 & 0 & 0 & -1 &  1 & 0 & 0 \\
       0 & 0 & 0 &  1 & 0 & 0 & 0 &  -1 & 0 \\
       0 & 0 &  -1 & 0 & 0 & 0 & 0 & 1 & 0 \\
       0 & 0 & 0 & 0 & 0 & 0 & 0 & 0 & 0 \\
       0 & 1 & 0 & 0 & 0 & 0 & -1 & 0 & 0 \\
       0 & -1 & 0 & 0 & 0 &  1 & 0 & 0 & 0 \\
       0 & 0 &  1 & -1 & 0 & 0 & 0 & 0 & 0 \\
       0 & 0 & 0 & 0 & 0 & 0 & 0 & 0  & 0
    \end{bmatrix},
\end{equation*}

\noindent which can also be written 

\begin{equation}
 h=\sum_{\gamma \in D_3}i (E^{\gamma(1)}_{\gamma(2)} 
\otimes E^{\gamma(2)}_{\gamma(3)} - E^{\gamma(2)}_{\gamma(3)} 
\otimes E^{\gamma(1)}_{\gamma(2)} )
\label{localham}
\end{equation}

\noindent where the elements $\gamma$ of $D_3$ are written as permutations of $\{1,2,3\}$. The above integrable system describes a one-dimensional lattice of anyons with $D(D_3)$ or $D(D_6)$ symmetry and local interactions given by (\ref{localham}).
The  two-site Hamiltonian may also be expressed in terms of spin-1 operators:

\begin{align*}
-ih = & - \sigma_+^2 \otimes (\sigma_z \sigma_- + \sigma_- \sigma_z) 
	+ (\sigma_z \sigma_- + \sigma_- \sigma_z) \otimes \sigma_+^2 \\
    & + \sigma_-^2 \otimes (\sigma_z \sigma_+ + \sigma_+ \sigma_z) 
	- (\sigma_z \sigma_+ + \sigma_+ \sigma_z) \otimes \sigma_-^2 \\
    & + \sigma_+ \sigma_z \otimes \sigma_z \sigma_+ - \sigma_z \sigma_+ \otimes \sigma_+ \sigma_z \\
    & + \sigma_- \sigma_z \otimes \sigma_z \sigma_- - \sigma_z \sigma_- \otimes \sigma_- \sigma_z 
\end{align*}

\noindent where $\sigma_{\pm} = \frac{1}{2} (\sigma_x \pm i\sigma_y)$.  This Hamiltonian is not the same as other known integrable spin-1 Hamiltonians \cite{kennedy,yung,ab}.

\section{Solutions of the Yang--Baxter equation associated with $D(A_4)$}

The same procedure can be applied to the symmetric and alternating groups.  The {\it symmetic group} $S_n$ is the group of permutations of $\{1, 2, ..., n\}$ where the operation is composition.  The subgroup of $S_n$ consisting of permutations which can be written as the product of an even number of transpositions is known as the {\it alternating group} and denoted $A_n$.  Now $S_3 \cong D_3$ and $A_3 \cong Z_3$, so we only consider $n \geq 4$.  In $A_n$ the only conjugacy class with only one element is $\{e\}$, which always gives rise to the trivial $R$-matrix $R= I \otimes I$.  Moreover, there are no conjugacy classes with two elements and only $A_4$ has a conjugacy class with three elements.
Therefore only $A_4$ can give rise to a three-dimensional irrep, and we can never obtain a two-dimensional irrep. 

Consider $A_4,$ using the convention $(12) \circ (13) = (132)$.  The relevent conjugacy class $C_k$ and the details required to construct the representations are:

\begin{align*}
&C_k = \{(12)(34), (13)(24), (14)(23)\},\quad g_k = (12)(34),\\
&Z_k = \{e, (12)(34), (13)(24), (14)(23)\}, \\
&\alpha_{(12)(34)} = e, \quad \alpha_{(13)(24)}=(132),\quad \alpha_{(14)(23)}=(123).
\end{align*}

This time $Z_k \cong \mathbb{Z}_2 \times \mathbb{Z}_2$, with the 4 one-dimensional irreps given by $(12)(34) = (-1)^a$, $(13)(24)=(-1)^b$, $a,b \in \{0,1\}$.  We obtain

$$R = {\rm diag} ((-1)^a, (-1)^{a+b}, (-1)^b, (-1)^b, (-1)^a, (-1)^{a+b}, (-1)^{a+b}, (-1)^b, (-1)^a).$$

\noindent Applying the permutation operator, we find

\begin{equation*}
\check{R} = \begin{bmatrix}
  (-1)^a & 0 & 0 & 0 & 0 & 0 & 0 & 0 & 0 \\
  0 & 0 & 0 & (-1)^b & 0 & 0 & 0 & 0 & 0 \\
  0 & 0 & 0 & 0 & 0 & 0 & (-1)^{a+b} & 0 & 0 \\
  0 & (-1)^{a+b} & 0 & 0 & 0 & 0 & 0 & 0 & 0 \\
  0 & 0 & 0 & 0 & (-1)^a & 0 & 0 & 0 & 0 \\
  0 & 0 & 0 & 0 & 0 & 0 & 0 & (-1)^b & 0 \\
  0 & 0 & (-1)^b & 0 & 0 & 0 & 0 & 0 & 0 \\
  0 & 0 & 0 & 0 & 0 & (-1)^{a+b} & 0 & 0 & 0 \\
  0 & 0 & 0 & 0 & 0 & 0 & 0 & 0  & (-1)^a
      \end{bmatrix},
\end{equation*}

\noindent which has eigenvalues $1,-1$ with multiplicities $6,3$ respectively when $a=0$, and eigenvalues $-1,i,-i$ 
each with multiplicity $3$ when $a=1$.  When $a=b=0$ this is the permutation matrix and gives rise to a representation of the Hecke algebra. Baxterisation then leads to the known $su(3)$ invariant solution. In the case $a=0, b=1$ $\check{R}$ is again a Hecke algebra representation, which can be Baxterised as $\check{R}(u) = I \otimes I + u \check{R}$.  This solution corresponds to the rational 15-vertex solution with a Reshetikhin twist \cite{reshetikhin}.  The last case is when $b=1$, in which case we can write without loss of generality:

\begin{equation*}
\check{R} = \begin{bmatrix}
  1 & 0 & 0 & 0 & 0 & 0 & 0 & 0 & 0 \\
  0 & 0 & 0 & 1 & 0 & 0 & 0 & 0 & 0 \\
  0 & 0 & 0 & 0 & 0 & 0 & -1& 0 & 0 \\
  0 & -1& 0 & 0 & 0 & 0 & 0 & 0 & 0 \\
  0 & 0 & 0 & 0 & 1 & 0 & 0 & 0 & 0 \\
  0 & 0 & 0 & 0 & 0 & 0 & 0 & 1 & 0 \\
  0 & 0 & 1 & 0 & 0 & 0 & 0 & 0 & 0 \\
  0 & 0 & 0 & 0 & 0 & -1 & 0 & 0 & 0 \\
  0 & 0 & 0 & 0 & 0 & 0 & 0 & 0  & 1
      \end{bmatrix}.
\end{equation*} 
We recognise this solution as belonging to the class of {\it non-standard} solutions of the 
Yang--Baxter equation. However, it is curious that the origin of the non-standard structure cannot be explained in terms of an underlying Lie superalgebra or colour Lie algebra structure \cite{mcanally, mb}, nor is it due to a Reshetikhin twist
\cite{reshetikhin}. 

As $\check{R}$ has three eigenvalues, we try the ansatz $\check{R}(x) = f(x) I \otimes I + g(x) \check{R} + h(x) \check{R}^{-1}$.  Applying the conjecture given in \cite{kauffman,cgx} we obtain 3 possible solutions.  Two of these, however, have $f(x)=0$ which is undesirable if we want the regularity property to hold. The third possible case can be shown to not satisfy the braiding Yang--Baxter equation
(\ref{mbybe}).  Hence we attempt to find another way to introduce a spectral parameter.

First we return to the original variables $u,v$ instead of $x,z$. 
Writing $a(u)= f(u) + g(u) + h(u)$ and $b(u)=g(u)-h(u)$, we note $\check{R}(u)$ is:

\begin{equation*}
\check{R}(u) = \begin{bmatrix}
 	a(u) & 0 & 0 & 0 & 0 & 0 & 0 & 0 & 0 \\
 	 0 & f(u) & 0 & b(u) & 0 & 0 & 0 & 0 & 0 \\
	 0 & 0 & f(u) & 0 & 0 & 0 & -b(u)& 0 & 0 \\
	 0 & -b(u) & 0 & f(u) & 0 & 0 & 0 & 0 & 0 \\
 	 0 & 0 & 0 & 0 & a(u) & 0 & 0 & 0 & 0 \\
 	 0 & 0 & 0 & 0 & 0 & f(u) & 0 & b(u) & 0 \\
 	 0 & 0 & b(u) & 0 & 0 & 0 & f(u) & 0 & 0 \\
 	 0 & 0 & 0 & 0 & 0 & -b(u) & 0 & f(u) & 0 \\
 	 0 & 0 & 0 & 0 & 0 & 0 & 0 & 0  & a(u)
      \end{bmatrix}.
\end{equation*}

\noindent Substituting $\check{R}(u)$ into the braiding Yang--Baxter equation (\ref{bybe}), 
we find $\check{R}(u)$ satisfies the Yang--Baxter equation if and only if the following conditions are met:

\begin{align}
&b(u+v) f(u) f(v) = f(u+v) [b(u) f(v) + b(v) f(u)], \label{b} \\
&a(u) b(u+v) f(v) = b(v) f(u) f(u+v) + a(u+v) b(u) f(v), \label{A1} \\
&a(u+v) f(u) f(v) = f(u+v) [a(u) a(v) + b(u) b(v)]. \label{A2}
\end{align}

First consider the case when $b(u) = 0$.  Then equations (\ref{b}, \ref{A1}) are automatically satisfied, and we need only consider equation \eqref{A2}.  The solution $\check{R}(u) \propto I \otimes I$ is uninteresting, so we instead choose $f(u)=1$ and $a(u)=e^u$, giving 

$$\check{R}(u) = {\rm diag}(e^u,1,1,1,e^u,1,1,1,e^u).$$
Observe that this solution has the following peculiar property:

\begin{equation}
\check{\R}=\lim_{u\rightarrow -\infty}\check{R}(u)\neq \check{R}. 
\label{weird}
\end{equation}

Next consider $b(u) \neq 0$.  We begin by choosing $f(u)=1$.  Then we see $b(u)=b_0 u$ is the only solution to \eqref{b}. We substitute these into equation (\ref{A1}) to obtain:

\begin{align*}
&(u+v)\;a(u) = v + u \;a(u+v),\\
&(u+v)\;a(v) = u + v \;a(u+v).
\end{align*}

Eliminating $a(u+v)$ we find:

\begin{align*}
&(u+v) \Bigl(\frac{a(u)}{u} - \frac{a(v)}{v} \Bigr) = \frac{v}{u} - \frac{u}{v} = \frac{v^2-u^2}{uv} \\
\Rightarrow \quad & \Bigl(\frac{a(u)}{u} - \frac{a(v)}{v} \Bigr) = \frac{v-u}{uv} = \frac{1}{u} - \frac{1}{v} \\
\Rightarrow \quad & \frac{a(u)-1}{u} = \frac{a(v)-1}{v} = c \\
\Rightarrow \quad & a(u) = 1 +cu.
\end{align*}
We find this satisfies \eqref{A2} provided $c = \pm ib_0$, so we have found a solution to the braiding Yang--Baxter equation (\ref{bybe}).  Note $b_0$ is just a scaling factor on $u$, so we can choose any non-zero complex number.  Choosing $b_0=i$ and $c=1$, we have

\begin{equation*}
\check{R}(u) = 
	\begin{bmatrix}
  1+u & 0 & 0 & 0 & 0 & 0 & 0 & 0 & 0 \\
  0 & 1 & 0 & iu & 0 & 0 & 0 & 0 & 0 \\
  0 & 0 & 1 & 0 & 0 & 0 & -iu& 0 & 0 \\
  0 & -iu & 0 & 1 & 0 & 0 & 0 & 0 & 0 \\
  0 & 0 & 0 & 0 & 1+u & 0 & 0 & 0 & 0 \\
  0 & 0 & 0 & 0 & 0 & 1 & 0 & iu & 0 \\
  0 & 0 & iu & 0 & 0 & 0 & 1 & 0 & 0 \\
  0 & 0 & 0 & 0 & 0 & -iu & 0 & 1 & 0 \\
  0 & 0 & 0 & 0 & 0 & 0 & 0 & 0  & 1+u
      \end{bmatrix}.
\end{equation*}  
The above solution again corresponds to the rational 15-vertex solution with a Reshetikhin twist \cite{reshetikhin}.
We also note that the property (\ref{weird}) also holds for this solution.

\section{Solutions of the Yang--Baxter equation associated with $D(S_4)$}

As with $D(A_n)$, the algebra $D(S_n)$ has no non-trivial $2$-dimensional irreps and the only non-trivial $3$-dimensional irrep occurs when $n=4$.  Then $C_k$, $g_k$, $Z_k$ and $\alpha_s, \, s \in C_k$ are given by:

\begin{align*}
&C_k = \{(12)(34), (13)(24), (14)(23)\},\quad g_k = (12)(34),\\
&Z_k = \{e,(12),(34),(12)(34), (13)(24), (14)(23), (1324),(1423)\}, \\
&\alpha_{(12)(34)} = e, \quad \alpha_{(13)(24)}=(14),\quad \alpha_{(14)(23)}=(13).
\end{align*}

Note that $Z_k \cong D_4$ with generators $\{ (12), (1324)\}$, so we know it has exactly $4$ one-dimensional reps given by $(12) = (-1)^a$, $(1324) = (-1)^b$, $a,b \in \{0,1\}$.  Following the same procedure as earlier, we obtain

\begin{equation*}
R = \sum_{g \in S_4} g \otimes g^* 
 = {\rm diag} ( 1, (-1)^{a+b}, (-1)^{a+b}, (-1)^{a+b},1,(-1)^{a+b},(-1)^{a+b},(-1)^{a+b}, 1).
\end{equation*}
Both these solutions arose in $D(A_4)$ and were discussed in the previous section.

\section{Conclusion}

Our results show that for certain constant solutions of the Yang--Baxter equation obtained by using representations of 
finite group doubles, it is possible to Baxterise them to yield solutions of the spectral parameter Yang--Baxter equation.
We have considered several examples where this is true and in particular we have found a new 21-vertex solution (\ref{newR}) from which we obtained an integrable model for a system of anyons with $D(D_3)$ or $D(D_6)$ symmetry.  
It is clearly of interest to determine if all constant solutions may be Baxterised. In contrast to the case of affine quantum algebras, where the spectral parameter has its origins in the loop representation, the origin of the spectral parameter for the above instances is unknown. 

In all our examples we have only looked for cases where the spectral parameter has the difference property. 
For the case of the generalised chiral Potts model in \cite{bkms} which does not have the difference property, an underlying finite group 
structure appears. This suggests that a Baxterisation ansatz without the assumption of the difference property may also be fruitful. Certainly more work is needed to fully realise the potential of finite group doubles in solving the 
Yang--Baxter equation with spectral parameter.  

\section*{Acknowledgements} This work was supported by the Australian Research Council.
We thank Mark Gould, Jabin Kirk and Liam Wagner for numerous helpful discussions on finite group doubles.
We also thank Carlos Mochon and Andrew Doherty for clarifying several aspects about properties of anyons, and Joost Slingerland for bringing relevant references to our attention. 

\section*{Appendix: Representations of $D(G)$ where $G=D_n$, $n$ odd} \label{irrepsodd}

The conjugacy classes $C_k$ of $G=D_n$, chosen representatives $g_k$, corresponding centraliser subgroups $Z_k$ and the 
elements $\alpha_s,\, \forall s \in C_k,$ are given below in Table \ref{odd}.

\begin{table}[ht]
\caption{$C_k, g_k$, $Z_k$ and $\alpha_s$ for $G=D_n$, $n$ odd.} \label{odd}
\centering
\begin{tabular}{|l|l|l|l|} \hline
$C_k$ & $g_k$ & $Z_k = Z(g_k)$ & $\alpha_s, \forall s \in C_k$ \\ \hline 
$\{e\}$ & $e$ & $D_n$ & $\alpha_e=e$\\
$\{\sigma^k, \sigma^{-k} \}, 1 \leq k \leq \frac{n-1}{2}$ & $\sigma^k$ & 
  $\{\sigma^i \,|\, 0 \leq i < n \}$ & $\alpha_{\sigma^k}=e, 
  \alpha_{\sigma^{-k}} = \tau$ \\
$\{ \sigma^i \tau \,|\, 0 \leq i < n \}$ & $\tau$ & $\{ e, \tau \}$ & 
  $\alpha_{\sigma^i \tau} = \sigma^{(\frac{n+1}{2})i}$\\ \hline
\end{tabular}
\end{table}

\noindent \underline{Representations induced by $C_k = \{e\}$}

\

The module elements are of the form $e \otimes v$ where $v \in V$, $V$ a $D_n$-module.  In representation terms, there are 2 one-dimensional irreps and $({n-1})/{2}$ two-dimensional irreps.  They are:

$$\sigma=1, \quad \tau = \pm 1, \quad g^* = \delta(g, e)$$

\noindent and 

\begin{equation*}
 \sigma = \begin{bmatrix}
                    \omega^k & 0  \\ 0 & \omega^{-k}
                 \end{bmatrix}, \quad
\tau = \begin{bmatrix}
                 0 & 1 \\ 1 & 0 
               \end{bmatrix}, \quad 
g^* = \delta(g, e) I_2
\end{equation*} 

\noindent where $1 \leq k \leq ({n-1})/{2}$.

\

\noindent \underline{Representations induced by $C_k = \{\sigma^k, \sigma^{-k}\}, \; 1 \leq k \leq ({n-1})/{2}$}

\

The module elements are of the form $e \otimes v, \tau \otimes v$ where $v \in V$, $V$ a module of the group algebra of $Z_k=\{ \sigma^i\, | \, 0 \leq i < n \}$.  There are $n$ such $A_k$ modules, with the corresponding representations $\pi_j$ given by $\pi_j (\sigma) = \omega^j,\, 0 \leq j < n$ where $\omega = \exp({{2\pi i}/{n}})$.  Thus there are $n(n-1)/{2}$ different irreps of $D(D_n)$ induced by these conjugacy classes, given by:

\begin{equation*}
 \sigma = \begin{bmatrix}
                    \omega^j & 0  \\ 0 & \omega^{-j}
                 \end{bmatrix}, \quad
\tau = \begin{bmatrix}
                 0 & 1 \\ 1 & 0 
               \end{bmatrix}, \quad 
(\sigma^k)^* = \begin{bmatrix}
                 1 & 0 \\ 0 & 0
               \end{bmatrix}, \quad
(\sigma^{-k})^* = \begin{bmatrix}
                 0 & 0 \\ 0 & 1
               \end{bmatrix}, \quad
g^* = 0 \text{ otherwise}
\end{equation*} 

\noindent where $0 \leq j < n, \, 1 \leq k \leq \frac{n-1}{2}$.

\

\noindent \underline{Representations induced by $C_k = \{\sigma^i \tau \,|\, 0 \leq i 
  < n \}$}

\

The module elements are of the form $\sigma^{j({n+1})/{2}} \otimes v,\: 0 \leq j < n$, 
where $v \in V$, $V$ a module of the group algebra of $Z_k = \{ e, \tau \}$.  
Hence there are two $n$-dimensional irreps of this form.  They are:

\begin{align*}
&\sigma = A \in M_{n \times n} \,\,\,\,\,\quad \text{ where } [A]_{ij} = \delta(i, j+2),
  \quad  \text{addition } mod \; n \\
&\tau =\pm B \in M_{n \times n} \quad \text{ where } [B]_{ij} = \delta(i+j, 2),
  \quad \text{addition } mod \; n \\
&(\sigma^i)^* = 0,\quad (\sigma^i \tau)^* = E^{i+1}_{i+1}, \quad 0 \leq i < n.
\end{align*}

\begin{example}
In $D(D_3)$, $\sigma$ and $\tau$ are as follows:
\begin{equation*}
\sigma = \begin{bmatrix}
           0 & 1 & 0 \\ 0 & 0 & 1 \\ 1 & 0 & 0 
	 \end{bmatrix}, \quad
\tau = \pm \begin{bmatrix}
             1 & 0 & 0 \\ 0 & 0 & 1 \\ 0 & 1 & 0
	   \end{bmatrix}.
\end{equation*}
\end{example}

Hence when $n$ is odd $D(D_n)$ has 2 one-dimensional irreps, 2 $n$-dimensional irreps and $({n^2-1})/{2}$ two-dimensional irreps, all of which are given above.  Note that the sum of the squares of the dimensions is $4n^2 = |D_n|^2 = |D (D_n)|$, as we expect \cite{gould93}.


\end{document}